\documentclass[12pt]{article}
\usepackage{latexsym}
\usepackage{amsmath}
\usepackage{amssymb}
\usepackage{amsfonts}
\numberwithin{equation}{section}
\newtheorem{theorem}{Theorem}[section]
\newtheorem{lemma}[theorem]{Lemma}

\newtheorem{claim}[theorem]{Claim}
\newtheorem{cor}[theorem]{Corollary}

\newcommand{\eproof}{{\mbox{\ }~\hfill
\mbox{\large $\Box$} \par \vskip 10pt}}

\newcommand{\ep}{\epsilon}

\newcommand{\R}{{\mathbf R}}

\newcommand{\pf}{\noindent{\bf Proof}}

\renewcommand{\d}{\partial}
\newcommand{\bd}{\bar\partial}

\title{Quantitative uniqueness estimates for second order elliptic equations with unbounded drift}
\author{Carlos Kenig\thanks{Department of Mathematics, University of Chicago, Chicago, IL 60637, USA. Email: cek@math.uchicago.edu. Supported in part by NSF Grant DMS-1265249 and DMS-0968472.}\qquad\quad Jenn-Nan Wang\thanks{Institute of Applied Mathematical Sciences, NCTS (Taipei), National Taiwan University,
Taipei 106, Taiwan. Email: jnwang@math.ntu.edu.tw. Supported in part by MOST Grants 102-2918-I-002-009 and 102-2115-M-002-009-MY3.}}

\date{}

\begin{document}
\maketitle

\begin{abstract}
In this paper we derive quantitative uniqueness estimates at infinity for solutions to an elliptic equation with unbounded drift in the plane.  More precisely, let $u$ be a real solution to $\Delta u+W\cdot\nabla u=0$ in $\R^2$, where $W$ is real vector and $\|W\|_{L^p(\R^2)}\le K$ for $2\le p<\infty$. Assume that $\|u\|_{L^{\infty}(\R^2)}\le C_0$ and satisfies certain a priori assumption at $0$. Then $u$ satisfies the following asymptotic estimates at $R\gg 1$
\[
\inf_{|z_0|=R}\sup_{|z-z_0|<1}|u(z)|\ge \exp(-C_1R^{1-2/p}\log R)\quad\text{if}\quad 2<p<\infty
\] 
and
\[
\inf_{|z_0|=R}\sup_{|z-z_0|<1}|u(z)|\ge R^{-C_2}\quad\text{if}\quad p=2,
\] 
where $C_1>0$ depends on $p, K, C_0$, while $C_2>0$ depends on $K, C_0$ . Using the scaling argument in \cite{bk05}, these quantitative estimates are easy consequences of estimates of the maximal vanishing order for solutions of the local problem. The estimate of the maximal vanishing order is a quantitative form of the strong unique continuation property. 
\end{abstract}

\section{Introduction}\label{intro}

In this work we consider the Schr\"odinger operator with an unbounded drift term
\begin{equation}\label{5-1}
\Delta u+W\cdot\nabla u=0\quad\text{in}\quad \R^2,
\end{equation}
where $W=(W_1,W_2)$ is a real vector-valued functions with $L^p$ bound for $2\le p<\infty$. Here we are interested in the lower bound of the decay rate for any nontrivial solution $u$. When $p=\infty$, the problem is related to Landis' conjecture \cite{kl88}. That is, let $u$ be a solution of \eqref{5-1} with $\|W\|_{L^{\infty}(\R^2)}\le 1$ and $\|u\|_{L^{\infty}(\R^2)}\le C_0$ and $|u(z)|\le\exp(-C|x|^{1+})$ for some $C>0$, then $u$ is trivial. If one applies a suitable Carleman estimate to \eqref{5-1} and a scaling devise in \cite{bk05}, the best exponent one can get is $2$, namely, under the same conditions stated above except $|u(z)|\le\exp(-C|x|^{2+})$, then $u$ is trivial (see \cite{da12}, \cite{lw13} for quantitative forms of this result). Moreover, in \cite{da12}, the author constructed a Meshkov type example showing that the exponent $2$ is in fact optimal for \emph{complex-valued} $W$ and $u$.

In a recent paper \cite{klw14}, the authors studied Landis' conjecture for second order elliptic equations in the plane in the real setting, including \eqref{5-1} with real-valued $W$ and $u$. It was proved in \cite{klw14} that if $u$ is a real-valued solution of \eqref{5-1} satisfying $|u(z)|\le\exp(C_0|z|)$, $|\nabla u(0)|=1$, and $\|W\|_{L^{\infty}(\R^2)}\le 1$, then
\begin{equation}\label{e0}
\inf_{|z_0|=R}\sup_{|z-z_0|<1}|u(z)|\ge \exp(-CR\log R)\quad\text{for}\quad R\gg 1
\end{equation}
where $C$ depends on $C_0$. 

In this paper, we would like to study estimates like \eqref{e0} for $2\le p<\infty$. For complex-valued $W$ satisfying
\begin{equation}\label{wd}
|W(z)|\le C\langle z\rangle^{-s},\quad s\ge 0,
\end{equation}
where $\langle z\rangle=\sqrt{1+|z|^2}$, the lower bound of the decay rate for $u$ is $\exp(-R^{2-2s}f(\log R))$ for $s<1/2$ and is $\exp(-R\tilde f(\log R))$ for $s\ge 1/2$, where $f(\log R)$ and $\tilde f(\log R)$ are functions of $\log R$ which grow slower than any positive power of $R$ (see \cite{da12}, \cite{lw13}). Here our assumption on $W$ will be an integral bound rather than a pointwise bound as in \eqref{wd}. Precisely, we prove that
\begin{theorem}\label{t1}
Let $u\in W^{2,p}_{loc}(\R^2)$ be a real solution of \eqref{5-1} with $|u(z)|\le C_0$ for some $C_0>0$ with $2\le p<\infty$.
\begin{description}
  \item[(i)] Assume that $2<p<\infty$,
  \begin{equation}\label{wl2}
  \|W\|_{L^p(\R^2)}\le \tilde K
  \end{equation}
  and $|\nabla u(0)|=1$. Then
  \[
  \inf_{|z_0|=R}\sup_{|z-z_0|<1}|u(z)|\ge \exp(-CR^{1-2/p}\log R)
  \]
  for $R\gg 1$, where $C$ depends on $p, \tilde K$, and $C_0$.
  
  \item[(ii)] For $p=2$, if
   \begin{equation}\label{L21}
  \|W\|_{L^2(\R^2)}\le K
  \end{equation}
  and
  \[
  1\le\int_{B_1}|\nabla u|^2,
  \]
  then
   \begin{equation}\label{2est}
  \inf_{|z_0|=R}\sup_{|z-z_0|<1}|u(z)|\ge R^{-C}
  \end{equation}
  for $R\gg 1$, where $C>0$ depends on $K, C_0$.
  \end{description}
Hereafter, we denote $B_r(a)$ the ball of radius $r$ centered at $a$. When $a=0$, we simply denote $B_r(a)=B_r$.
\end{theorem}

Using the scaling argument in \cite{bk05}, Theorem~\ref{t1} is an easy consequence of the estimate of the maximal vanishing order of the solution $v$ to
\begin{equation}\label{local}
\Delta v+A\cdot\nabla v=0\quad\text{in}\quad B_8
\end{equation}
with 
\begin{equation}\label{am}
\|A\|_{L^p(B_8)}\le K. 
\end{equation}
It suffices to take $K\ge 1$. The proof of the maximal vanishing order of $v$ relies on a nice reduction of \eqref{local} to a $\bd$ equation. Having the $\bd$ equation, we then derive the vanishing order by using Hadamard's three circle theorem. The case $p=2$ needs special attention due to the fact that the Cauchy transform fails to be a bounded map from $L^2(B_8)$ to $L^{\infty}(B_8)$. 

The estimate of the maximal vanishing order of $v$ provides us a quantitative form of the strong unique continuation property (SUCP) for \eqref{local}. Note that $A\in L^2$ is a scale invariant drift in $\R^2$ in the sense that if $v(x)$ solves \eqref{local}, then $v_r(x):=v(rx)$ satisfies $\Delta v_r+A_r\nabla v_r=0$ in $B_{8/r}$ with $A_r(x)=rA(rx)$ and
\[
\|A\|_{L^2(B_8)}=\|A_r\|_{L^2(B_{8/r})}.
\]
It is clear that $v(z)=\exp(-|z|^{-\ep})$ for $\ep>0$ is an easy counterexample of SUCP for $A\in L^p$ with $p<2$. For the dimension $n\ge 3$, Kim \cite{ki89} proved that SUCP holds for \eqref{local} when $A\in L_{loc}^p$ with $p=(3n-2)/2$ and Wolff \cite{wo90} improved the exponent to $p=\max\{n,(3n-4)/2\}$.  On the other hand, if $n\ge 5$, counterexamples to the SUCP with $A\in L_{loc}^n$ were given by Wolff in \cite{wo94} (or see \cite{wo93}). Counterexamples of the unique continuation property (UCP) for \eqref{local} with $A\in L^p$, $p<2$, or $A\in L^2_{weak}$, weak $L^2$ space, were constructed by Mandache \cite{ma02} and Koch-Tataru \cite{kt02}, respectively. We also would like to mention that a counterexample of UCP for the Schr\"odinger operator $\Delta u+Vu=0$ with $V\in L^1$ was constructed by Kenig and  Nadirashvili \cite{kn00} for dimension $n\ge 2$. For $n=2$ and $A\in L^2$, it seems likely that a variant of the Carleman estimate proved in Kim's thesis for $n\ge 3$ \cite[Theorem~3]{ki89} is available for $n=2$ and the SUCP will follow from it (see the remark in \cite[Page 156]{wo90}). Here we provide an explicit proof of the SUCP for \eqref{local} in two dimensions, where $A\in L_{loc}^2$ is a real-valued vector. Using the same method, we also study the SUCP for 
\begin{equation}\label{diw}
\Delta v+\nabla\cdot(Av)=0, 
\end{equation}
where $A$ is a real-valued vector with bounded $L_{loc}^2$ norm. 

The structure of the paper is as follows. In Section~2, we consider the case where $2<p<\infty$. The case of $p=2$ is treated in Section~3. We study the SUCP for \eqref{diw} in Section~4. Throughout the paper, $C$ stands for an absolute constant whose dependence will be specified if necessary. Its value may vary from line to line. 

\section{The case of $2<p<\infty$}\label{sec2}

We consider any solution $v\in W_{loc}^{2,p}(B_8)$ to the equation \eqref{local} with $A=(A_1,A_2)$ satisfying \eqref{am}. Denote $g=v_x-iv_y$. It is easy to see that
\begin{equation}\label{5-2}
\bar\d g=\frac 12\Delta v=-\frac 12(A_1\d_xv+A_2\d_yv)=-\frac 14(A_1+iA_2)g-\frac 14(A_1-iA_2)\bar g,
\end{equation}
where in the last step we used that $v$ is real. As usual, we denote $\bd=(\d_x+i\d_y)/2$. Let us define
\[
\alpha(z)=\left\{
\begin{aligned}
&-\frac 14(A_1+iA_2)-\frac 14(A_1-iA_2)\frac{\bar g}{g}\quad&\text{if}\quad g\ne 0,\\
&0\quad&\text{if}\quad g=0,
\end{aligned}
\right.
\]
then \eqref{5-2} can be written as
\begin{equation}\label{5-3}
\bd g=\alpha g\quad\text{in}\quad B_8.
\end{equation}
Therefore, any solution of \eqref{5-3} is represented by
\begin{equation}\label{gwh}
g=\exp(w)h\quad\text{in}\quad B_8,
\end{equation}
where $h$ is holomorphic in $B_8$ and 
\begin{equation}\label{w}
w(z)=-\frac{1}{\pi}\int_{B_8}\frac{\alpha(\xi)}{\xi-z}d\xi,
\end{equation}
i.e., $w$ is the Cauchy transform of $\alpha$.

From \eqref{am} and the definition of $\alpha$, we have that
\[
\|\alpha\|_{L^p(B_8)}\le K
\]
with $2<p$. In view of the mapping properties of the Cauchy transform (see for example \cite{ve62}), we see that
\begin{equation}\label{wz}
|w(z)|\le CK\quad\text{for}\quad z\in B_8,
\end{equation}
where $C$ depends on $p$. Since $h$ is holomorphic in $B_8$, Hadamard's three circle theorem implies
\begin{equation*}
\|h\|_{L^{\infty}(B_{1})} \le \|h\|_{L^{\infty}(B_{r/4})}^{\theta}\|h\|_{L^{\infty}(B_{6})}^{1-\theta},
\end{equation*}
where we choose $r/4<1$ and
\[
\theta=\frac{\log(6)}{\log({24}/{r})}.
\]
Standard interior estimates imply that
\begin{equation}\label{3l2}
\|h\|_{L^{\infty}(B_{1})} \le C(r^{-1}\|h\|_{L^{2}(B_{r/2})})^{\theta}\|h\|_{L^{2}(B_{7})}^{1-\theta}.
\end{equation}
On the other hand, it is not hard to prove that $v$ of \eqref{local} satisfies the following Caccioppoli's inequality
\begin{equation}\label{caccio}
\int_{B_r}|\nabla v|^2\le\frac{C\|A\|_{L^p(B_8)}^2}{(\rho-r)^2}\|v\|^2_{L^{\infty}(B_\rho)}\le\frac{CK^2}{(\rho-r)^2}\|v\|^2_{L^{\infty}(B_\rho)},\quad 0<r<\rho<8,
\end{equation}
where $C$ depends on $p$. The derivation of \eqref{caccio} follows from the standard procedure using a cutoff function. We omit the details here. Combining \eqref{gwh}, \eqref{3l2} and \eqref{caccio}, we have that
\begin{equation}\label{est0}
\exp(-CK)\|\nabla v\|_{L^{\infty}(B_1)}\le C(r^{-1}\exp(CK)\|v\|_{L^{\infty}(B_{r})})^{\theta}\|v\|_{L^{\infty}(B_{8})}^{1-\theta}.
\end{equation}
Based on \eqref{est0}, we immediately prove 
\begin{theorem}\label{tt1}
Let $v\in W_{loc}^{2,p}(B_8)$ be a real solution of \eqref{local} with $A$ satisfying \eqref{am}. Assume that $v$ satisfies $|v(z)|\le C_0$ for all $z\in B_8$ and $\sup_{B_1}|\nabla v(z)|\ge 1$. Then 
\begin{equation}\label{rk}
\|v\|_{L^{\infty}(B_r)}\ge r^{C_1+C_2K},
\end{equation}
where $C_1$ depends on $C_0$ and $C_2$ depends on $p$.
\end{theorem}

From Theorem~\ref{tt1}, we can easily derive the following quantitative uniqueness estimate, which is (i) of Theorem~\ref{t1}.
\begin{cor}\label{cc1}
Let $u\in W_{loc}^{2,p}(\R^2)$ be a real solution of \eqref{5-1} with $|u(z)|\le C_0$ and $|\nabla u(0)|=1$. Assume that
\[
\|W\|_{L^p(\R^2)}\le \tilde K.
\]
Then 
\begin{equation}\label{p1}
\inf_{|z_0|=R}\sup_{|z-z_0|<1}|u(z)|\ge \exp(-CR^{1-2/p}\log R)
\end{equation}
for $R\gg 1$, where $C$ depends on $p, \tilde K$ and $C_0$.
\end{cor}
\pf \, We use the scaling argument in \cite{bk05}.  Precisely, let $|z_0|=R$ with $R\gg 1$, and define $u_R(z)=u(R(z+z_0/R))$. Then $u_R$ satisfies
\[
\Delta u_R+W_R\cdot\nabla u_R=0\quad\text{in}\quad B_8,
\]
where $W_R(z)=RW(R(z+z_0/R))$. It is clear that
\[
\left(\int_{B_8}|W_R|^p\right)^{1/p}\le R^{1-2/p}\left(\int_{\R^2}|W|^p\right)^{1/p}\le \tilde KR^{1-2/p}.
\]
Also, we observe that
\[
|\nabla u_R(-{z_0}/{R})|=R|\nabla u(0)|=R>1.
\]
Taking $K=\tilde KR^{1-2/p}$ and $r=R^{-1}$, estimate \eqref{rk} yields \eqref{p1}.\eproof

\section{The case of $p=2$}\label{sec3}

Likewise, we consider the local problem \eqref{local}. Here we assume that 
\begin{equation}\label{2am}
\|A\|_{L^2(B_8)}\le K.
\end{equation}
We first establish an estimate of the maximal vanishing order of $v$ to \eqref{local} under the assumption \eqref{2am}. 
\begin{theorem}\label{tt12}
Let $v\in W^{2,2}_{loc}(B_8)$ be a real solution of \eqref{local} with $A$ satisfying \eqref{2am}. Assume that $v$ satisfies $|v(z)|\le C_0$ for all $z\in B_8$ and 
\[
\|\nabla v\|_{L^2(B_{6/5})}\ge 1. 
\]
Then for $r$ small
\begin{equation}\label{rkk}
\|v\|_{L^{\infty}(B_r)}\ge r^{C_1+C_2K^2},
\end{equation}
where $C_1$ depends on $C_0$ and $C_2$ is an absolute constant. 
\end{theorem}

The proof of Theorem~\ref{tt12} is more involved. Note that the formula \eqref{5-3} remains valid, i.e.,
\begin{equation}\label{0701}
\bd g=\alpha g\quad\text{in}\quad B_8,
\end{equation}
and
\[
\|\alpha\|_{L^2(B_8)}\le K.
\]
Likewise, let
\begin{equation*}
w(z)=\frac{1}{\pi}\int_{B_8}\frac{\alpha(\xi)}{\xi-z}d\xi,
\end{equation*}
then any solution of \eqref{0701} is represented by
\[
g(z)=\exp(-w(z))h(z)\quad\text{for}\quad z\in B_8
\]
where $h$ is holomorphic in $B_8$. It is not hard to see that
\[
\|w\|_{W^{1,2}(B_8)}=\|w\|_{L^2(B_8)}+\|\nabla w\|_{L^2(B_8)}\le CK.
\]

In the sequel, we need to estimate $\int_{B_r}\exp(2|w|)$ for $r\le 2$. For this end, we recall the following Trudinger's Sobolev embedding theorem in the plane \cite{st72}, \cite{tr67}. Assume that $f\in W^{1,2}(B_1)$ and $\|f\|_{W^{1,2}(B_1)}\le 1$, then there exist two absolute constants $\tilde\alpha_\ast$ and $\tilde C_\ast$ such that
\[
\int_{B_1}\exp(\tilde\alpha_\ast f^2)\le \tilde C_\ast.
\]
By Poincar\'e's inequality, we immediately obtain that
\begin{cor}\label{cor02}
If $f\in W^{1,2}(B_1)$, $\int_{B_1}f=0$, and $\|\nabla f\|_{L^2(B_1)}\le 1$, then there exist $\alpha_\ast$ and $C_\ast$ such that
\[
\int_{B_1}\exp(\alpha_\ast f^2)\le C_\ast.
\]
\end{cor}

Our task now is to prove 
\begin{lemma}\label{lemma0701}
For $q>0$ and $0<r\le 2$, we have that
\begin{equation}\label{qqest}
\frac{1}{|B_r|}\int_{B_r}\exp(q|w|)\le Cr^{-qCK}\exp(qCK+q^2CK^2).
\end{equation}
\end{lemma}
\pf\, By a scaling argument, we can deduce from Corollary~\ref{cor02} that if $f\in W^{1,2}(B_r)$, $\int_{B_r}f=0$, and $\|\nabla f\|_{L^2(B_r)}\le 1$, then
\begin{equation}\label{q0701}
\frac{1}{r^2}\int_{B_r}\exp(\alpha_\ast f^2)\le C_\ast.
\end{equation}
To verify \eqref{q0701}, we define $f_r(x)=f(rx)$ for $x\in B_1$ and observe that
\[
\int_{B_1}|\nabla f_r|^2=\int_{B_r}|\nabla f|^2.
\]
Then \eqref{q0701} follows directly from Corollary~\ref{cor02}. 

Let us define $w_r(x)=w(x)-\overline{w}_r$, where $\overline{w}_r=\frac{1}{|B_r|}\int_{B_r}w$. We first consider the case when $\|\nabla w_r\|_{L^2(B_r)}>0$.  We can write
\begin{equation}\label{07011}
\int_{B_r}\exp(q|w_r|)=\int_{B_r}\exp\left(q\|\nabla w_r\|_{L^2(B_r)}\cdot\left|\frac{w_r}{\|\nabla w_r\|_{L^2(B_r)}}\right|\right)=\int_{B_r}\exp(a|f|),
\end{equation}
where
\[
a=q\|\nabla w_r\|_{L^2(B_r)}\quad\text{and}\quad f=\frac{w_r}{\|\nabla w_r\|_{L^2(B_r)}}.
\]
Note that $\|\nabla w_r\|_{L^2(B_r)}\le CK$. It is helpful to study the function $e^{ax}$ for $x>0$. We first consider the case when $ax\le\alpha_\ast x^2$, i.e., $x\ge a/\alpha\ast$. In this case, it is trivial that $e^{ax}\le e^{\alpha_\ast x^2}$. In the case when $x\le a/\alpha_\ast$, we have $e^{ax}\le e^{a^2/\alpha_\ast}$. Consequently, we obtain that
\[
e^{ax}\le e^{\alpha_\ast x^2}+e^{a^2/\alpha_\ast},\quad x>0.
\] 
Therefore, it follows from \eqref{q0701} and \eqref{07011} that 
\begin{equation}\label{06241}
\begin{aligned}
\int_{B_{r}}\exp(q|w_r|)&\le\int_{B_r}\exp(\alpha_\ast|f|^2)+\int_{B_r}\exp(a^2/\alpha_\ast)\\
&\le (C_\ast+\exp(a^2/\alpha_\ast))r^2\le Cr^2\exp(q^2CK^2).
\end{aligned}
\end{equation}

Next we want to estimate $|\overline{w}_r|$. 
\begin{claim}\label{claim1}
\[
|\overline{w}_{r}|\le CK\log(1/r)+CK.
\]
\end{claim}
\pf\, Note that
\[
\begin{aligned}
|\overline{w}_r-\overline{w}_{2r}|&=\left|\frac{1}{|B_r|}\int_{B_r}w-\overline{w}_{2r}\right|\le\frac{1}{|B_r|}\int_{B_r}|w-\overline{w}_{2r}|\\
&\le\frac{C}{|B_{2r}|}\int_{B_{2r}}|w-\overline{w}_{2r}|\le C\left(\frac{1}{|B_{2r}|}\int_{B_{2r}}|w-\overline{w}_{2r}|^2\right)^{1/2}\\
&\le C\left(\int_{B_{2r}}|\nabla w|^2\right)^{1/2}\le CK.
\end{aligned}
\]
It is clear that
\begin{equation}\label{0611}
|\overline{w}_{r}|\le |\overline{w}_{r}-\overline{w}_{2r}|+ |\overline{w}_{2r}-\overline{w}_{4r}|+\cdots+|\overline{w}_{2^kr}|.
\end{equation}
We now choose 
\[
k=\lfloor\frac{1}{\log 2}\log(\frac 1r)\rfloor+1\le C\log(\frac 1r),
\]
where $\lfloor\cdot\rfloor$ is the floor function. With the choice of $k$, we can see that
\[
1\le 2^kr\le 2.
\]
Each term of \eqref{0611} is bounded by $CK$. The claim follows immediately.\eproof

It is clear that Claim~\ref{claim1} implies
\begin{equation}\label{wbar}
\exp(q|\overline{w}_r|)\le \exp(qCK)r^{-qCK}.
\end{equation}
Combining \eqref{06241} and \eqref{wbar} yields
\[
\int_{B_r}\exp(q|w|)\le\int_{B_r}\exp(q|w-\overline{w}_r|)\exp(q|\overline{w}_r|)\le Cr^{2-qCK}\exp(qCK+q^2CK^2).
\]
Now if $\|\nabla w_r\|_{L^2(B_r)}=0$, then $w(x)\equiv\overline{w}_r$  in $B_r$. Hence, we have
\[
\int_{B_r}\exp(q|w|)=\int_{B_r}\exp(q|\overline{w}_r|)\le Cr^{2-qCK}\exp(qCK).
\]
The derivation of \eqref{qqest} is now completed.\eproof

As above, we will apply Hadamard's three circle theorem to $h=\exp(w)g$ with $r_2=6/5$, $r_3=2$, and $r_1=r/4<6/5$, i.e.,
\begin{equation}\label{e120}
\|\exp(w)g\|_{L^{\infty}(B_{r_2})} \le \|\exp(w)g\|_{L^{\infty}(B_{r_1})}^{\theta}\|\exp(w)g\|_{L^{\infty}(B_{r_3})}^{1-\theta},
\end{equation}
where
\begin{equation}\label{th}
\theta=\frac{\log(10/6)}{\log(8/r)}.
\end{equation}
We will estimate the terms on both sides of \eqref{e120}. We begin with the terms on the right hand side. Note that $v$ here also satisfies Caccioppoli's estimate \eqref{caccio} for $p=2$. On the other hand, using the Poisson kernel of the unit disc, it is easy to see that for any holomorphic function $h$
\[
\|h\|_{L^\infty(B_{r/2})}\le C\frac{1}{|B_r|}\int_{B_r}|h|.
\]
Putting all estimates together and in view of $g=v_x-iv_y$, we have that 
\begin{equation}\label{rest1}
\begin{aligned}
\|\exp(w)g\|_{L^{\infty}(B_{r/4})}&=\|h\|_{L^{\infty}(B_{r/4})}\le \frac{C}{|B_{r/2}|}\int_{B_{r/2}}|\exp(w)g|\\
&\le C\left(\frac{1}{|B_{r/2}|}\int_{B_{r/2}}\exp(2|w|)\right)^{\frac{1}{2}}\left(\frac{1}{|B_{r/2}|}\int_{B_{r/2}}|\nabla v|^{{2}}\right)^{\frac{1}{2}}\\
&\le Cr^{-CK}\exp(CK^2)\|\nabla v\|_{L^{2}(B_{r/2})}\\
&\le CKr^{-CK}\exp(CK^2)\|v\|_{L^{\infty}(B_{r})}\le C^{CK^2}r^{-CK}\|v\|_{L^\infty(B_r)},
\end{aligned}
\end{equation}
where we used \eqref{qqest} with $q=2$ in the third inequality and Caccioppoli's estimate in the fourth inequality. Using \eqref{rest1} on the right hand side of \eqref{e120} gives
\begin{equation}\label{0702}
\begin{aligned}
\|\exp(w)g\|_{L^{\infty}(B_{r/4})}^{\theta}\|\exp(w)g\|_{L^{\infty}(B_{2})}^{1-\theta}&\le (C^{CK^2}r^{-CK}\|v\|_{L^\infty(B_r)})^{\theta}(C^{CK^2}8^{-CK}\|v\|_{L^\infty(B_8)})^{1-\theta}\\
&\le C_0C^{CK^2}(C^{CK^2}r^{-CK}\|v\|_{L^\infty(B_r)})^{\theta}.
\end{aligned}
\end{equation}

We now turn to the estimate of $\|\exp(w)g\|_{L^{\infty}(B_{r_2})}=\|\exp(w)g\|_{L^{\infty}(B_{6/5})}$ on the left side of \eqref{e120}. From \eqref{qqest} with $q=4$ and $r=6/5$, it is readily seen that
\begin{equation}\label{353}
\begin{aligned}
&1\le\|\nabla v\|_{L^{2}(B_{6/5})}=\|g\|_{L^{2}(B_{6/5})}=\|\exp(-w)h\|_{L^{2}(B_{6/5})}\\
&\le \|\exp(|w|)\|_{L^{4}(B_{6/5})}\|h\|_{L^{4}(B_{6/5})}\le C^{CK^2}\|h\|_{L^{\infty}(B_{6/5})}.
\end{aligned}
\end{equation}
Combining \eqref{0702}, \eqref{353} and the form of $\theta$ (see \eqref{th}), we immediately arrive at the estimate \eqref{rkk}. The proof of Theorem~\ref{tt12} is completed.

Now we can put everything together to prove (ii) of Theorem~\ref{t1}. 

\medskip\noindent
\emph{Proof of (ii) of Theorem~\ref{t1}}. Let $|z_0|=R\gg 1$ and $v(z)=u(R(z+z_0/R))$. Then $v$ solves \eqref{local} and with $A(z)=RW(R(z+z_0/R))$. Note that
\[
\|A\|_{L^2(B_8)}\le K
\]
since $\|W\|_{L^2(\R^2)}\le K$. The boundedness assumption on $u$ implies $\|v\|_{L^{\infty}(B_8)}\le C_0$. On the other hand, we can see that for $\tilde z_0=-z_0/R$ ($|\tilde z_0|=1$)
\begin{equation*}
1\le\|\nabla u\|_{L^2(B_1)}=\|\nabla v\|_{L^2(B_{1/R}(\tilde z_0))}\le\|\nabla v\|_{L^2(B_{6/5})}\end{equation*}
provided $R$ is large. Therefore, letting $r=1/R$ in \eqref{rkk}, we obtain that
\begin{equation*}
\|u\|_{L^{\infty}(B_1(z_0))}=\|v\|_{L^{\infty}(B_{r}(0))}\ge r^{C_1}=R^{-C_1},
\end{equation*}
where $C_1>0$ depends on $C_0$ and $K$. \eproof

Note that $v(z)-v(0)$ is also a solution of \eqref{local}. Thus the estimate of vanishing order \eqref{rkk} remains valid for $v(z)-v(0)$. Consequently, we obtain the following (SUCP) result.
\begin{cor}\label{sucp}
Assume that $\Omega$ is an open connected domain of $\R^2$. Let $v\in W_{loc}^{2,2}(\Omega)$ be any solution of 
\[
\Delta v+A\cdot\nabla v=0\quad\text{in}\quad\Omega,
\]
with real-valued drift $A\in L^2(\Omega)$, then $v$ satisfies (SUCP), namely, if for some $z_0\in\Omega$
\begin{equation*}
|v(z)-v(z_0)|={\mathcal O}(|z-z_0|^N)\;\;\text{for all}\;\; N\in{\mathbb N},\;\;\text{as}\;\; |z-z_0|\to 0,
\end{equation*}
i.e., if for $N\in{\mathbb N}$, there exist $C_N>0$ and $r_N>0$ such that
\[
|v(z)-v(z_0)|\le C_N|z-z_0|^N\quad\forall\;\; |z-z_0|<r_N,
\]
then $v(z)\equiv v(z_0)$ for all $z\in\Omega$.  
\end{cor}
\pf\; It suffices to consider a real solution $v$. First assume that $z_0=0$ and $B_8\subset\Omega$. We can always assume this by translation and scaling. Note that $\|A\|_{L^2(B_8)}$ is finite. If $v(z)\not\equiv v(0)$ in $B_{6/5}$, then $\|\nabla v\|_{L^2(B_{6/5})}\ge C$
for some $C>0$. The estimate \eqref{rkk} implies that $v(z)-v(0)$ cannot vanish at $0$ to infinite order. Therefore, we must have $v(z)=v(0)$ for all $z\in B_{6/5}$. A chain of balls argument then finishes the proof.
\eproof

\section{SUCP for an equation of divergence form}

In this section, we would like to prove the SUCP for solutions of
\begin{equation}\label{divw0}
\Delta v+\nabla\cdot(Av)=0\quad\text{in}\quad\Omega,
\end{equation}
where $\Omega\subset\R^2$ is an open connected domain and $A=(A_1,A_2)$ is a real-valued vector satisfying 
\begin{equation}\label{a1}
\|A\|_{L^2(\Omega)}\le C_0.
\end{equation}
In other words, we will show that
\begin{theorem}\label{sucp2}
Let $v\in W_{loc}^{1,2}(\Omega)$ be any solution of \eqref{divw0}. Let $z_0\in\Omega$ and
\[
|v(z)|={\cal O}(|z-z_0|^N)\quad\text{as}\quad|z-z_0|\to 0
\]
for all $N>0$, then $v\equiv 0$ in $\Omega$.
\end{theorem}
\pf\; As before, it suffices to consider a real solution $v$. We first assume $z_0=0$, $B_8\subset\Omega$ and consider
\begin{equation}\label{divw}
\Delta v+\nabla\cdot(Av)=0\quad\text{in}\quad B_8.
\end{equation}
Since \eqref{divw} is of divergence form, there exists $\tilde v$ with $\tilde v(0)=0$ such that
\begin{equation}\label{vv}
\begin{cases}
\d_y\tilde v=\d_xv+A_1v,\\
-\d_x\tilde v=\d_yv+A_2v.
\end{cases}
\end{equation}
Let $f=v+i\tilde v$, then $f$ satisfies
\begin{equation}\label{ff}
\bd f=\frac 12(A_1+iA_2)v=\frac 14(A_1+iA_2)(f+\bar f)=\alpha f,
\end{equation}
where
\[
\alpha=\left\{
\begin{aligned}
&\frac 14(A_1+iA_2)(1+\frac{\bar f}{f})\quad\text{if}\quad f\ne 0,\\
&0\quad\text{if}\quad f=0.
\end{aligned}\right.
\]
It follows from \eqref{a1} that
\begin{equation}\label{a2}
\|\alpha\|_{L^2(B_8)}\le C_0.
\end{equation}
Any solution of \eqref{ff} in $B_8$ is written as $f=\exp(-w)h$, where $h$ is holomorphic in $B_8$ and 
\[
w(z)=\frac{1}{\pi}\int_{B_8}\frac{\alpha(\xi)}{\xi-z}d\xi.
\]
As before, we have that
\[
\|w\|_{W^{1,2}(B_8)}\le C,
\]
where $C$ depends on $C_0$.

Applying Hadamard's three circle theorem to $h=\exp(w)f$ with $r_1=r/4<1$, $r_2=1$, $r_3=2$, we have that
\begin{equation}\label{e620}
\|\exp(w)f\|_{L^{\infty}(B_1)} \le \|\exp(w)f\|_{L^{\infty}(B_{r/4})}^{\theta}\|\exp(w)f\|_{L^{\infty}(B_{2})}^{1-\theta},
\end{equation}
where
\begin{equation*}
\theta=\theta(r)=\frac{\log 2}{\log(8/r)}.
\end{equation*}
As in the estimate \eqref{353}, we can see that
\begin{equation}\label{3533}
\|v\|_{L^{2}(B_{1})}\le\|f\|_{L^{2}(B_{1})}=\|\exp(-w)h\|_{L^{2}(B_{1})}\le \|\exp(|w|)\|_{L^{4}(B_{1})}\|h\|_{L^{4}(B_{1})}\le C\|h\|_{L^{\infty}(B_{1})}.
\end{equation}
This estimate will give us a lower bound on the right hand side of \eqref{e620}.

It is not hard to prove that a Caccioppoli's type inequality holds for the solution $v$ of \eqref{divw}, i.e., for $r<\rho<8$, we have
\begin{equation}\label{caccio2}
\int_{B_r}|\nabla v|^2\le\frac{C}{(\rho-r)^2}\|v\|^2_{L^{\infty}(B_\rho)}.
\end{equation}
As in the derivation of \eqref{rest1}, we can obtain that
\begin{equation}\label{rest27}
\begin{aligned}
\|\exp(w)f\|_{L^{\infty}(B_{r/4})}&=\|h\|_{L^{\infty}(B_{r/4})}\le \frac{C}{|B_{r/2}|}\int_{B_{r/2}}|\exp(w)f|\\
&\le C\left(\int_{B_{r/2}}\exp(2|w|)\right)^{\frac{1}{2}}\left(\int_{B_{r/2}}|f|^{2}\right)^{\frac{1}{2}}\\
&\le Cr^{-C}\left(\int_{B_{r/2}}|f|^{2}\right)^{\frac{1}{2}}\le Cr^{-C}(\|v\|_{L^2(B_{r/2})}+\|\tilde v\|_{L^2(B_{r/2})}),
\end{aligned}
\end{equation}
where $0<r<8$. We now need to estimate $\|\tilde v\|_{L^2(B_{r/2})}$ in \eqref{rest27}. To this end, we can use \eqref{vv} and \eqref{caccio2} to compute

\noindent
\begin{equation}\label{0527}
\begin{aligned}
\int_{B_{r/2}}|\tilde v(x)|^2&=\int_{B_{r/2}}|\tilde v(x)-\tilde v(0)|^2=\int_{B_{r/2}}|\int_0^1\nabla\tilde v(tx)\cdot xdt|^2dx\\
&\le (r/2)^2\int_{B_{r/2}}\int_0^1|\nabla\tilde v(tx)|^2dtdx\\
&\le Cr^3\int_0^{r/2}\left\{\frac{1}{|B_s|}\int_{B_s}|\nabla\tilde v(y)|^2dy\right\}ds\\
&\le Cr^3\int_0^{r/2}\left\{\frac{1}{|B_s|}\int_{B_s}(|\nabla v(y)|^2+|Av|^2)dy\right\}ds\\
&\le Cr^3\int_0^{r/2}\left\{\frac{\|v\|^2_{L^{\infty}(B_{2s})}}{s^2|B_s|}+\frac{\|v\|^2_{L^{\infty}(B_{s})}}{|B_s|}\int_{B_s}|A|^2dy\right\}ds\\
&\le Cr^3\int_0^{r/2}\left\{\frac{\|v\|^2_{L^{\infty}(B_{2s})}}{s^2|B_s|}+\frac{\|v\|^2_{L^{\infty}(B_{s})}}{|B_s|}\right\}ds.
\end{aligned}
\end{equation}

The assumption that $v$ vanishes at $0$ to infinite order implies that there exist $C_4>0$ and $r_4<8$ such that
\[
|v(z)|\le C_4|z|^4,\quad\forall\;\;|z|<r_4.
\]
The estimate \eqref{0527} gives us
\begin{equation}\label{05277}
\int_{B_4}|\tilde v|^2\le C\left(\int_0^{r_4/2}+\int_{r_4/2}^8\right)\left\{\frac{\|v\|^2_{L^{\infty}(B_{2s})}}{s^2|B_s|}+\frac{\|v\|^2_{L^{\infty}(B_{s})}}{|B_s|}\right\}ds\le C.
\end{equation}
Combining \eqref{rest27} and \eqref{05277} yields
\begin{equation}\label{ff2}
\|\exp(w)f\|^{1-\theta}_{L^{\infty}(B_{2})}\le C'
\end{equation}
for all $0<\theta<1$, where $C'>0$. Now if we assume that
\begin{equation}\label{b1}
\|v\|_{L^{2}(B_{1})}\ge e^{-k}
\end{equation}
for some $k>0$, then we obtain from \eqref{e620}, \eqref{3533}, and \eqref{ff2} that
\begin{equation*}
\tilde Cr^{\tilde C k} \le \|\exp(w)f\|_{L^{\infty}(B_{r/4})},
\end{equation*}
where $\tilde C$ depends on $C'$. However, using the fact that $v$ vanishes at $0$ to infinite order, \eqref{rest27}, \eqref{0527}, we have that there exist $N_0>\tilde Ck$ and $r_{N_0}$ so that
\[
\|\exp(w)f\|_{L^{\infty}(B_{r/4})}\le C_{N_0}r^{N_0}
\]
for all $r<r_{N_0}$. This leads to a contradiction. In other words, we must have $\|v\|_{L^{2}(B_{1})}< e^{-k}$ for all $k>0$ and hence $v\equiv 0$ in $B_1$. 

Now we consider the general case, i.e., $v$ vanishes at some $z_0\in\Omega$ to infinite order. We choose a $r_0$ satisfying $B_{8r_0}(z_0)\subset\Omega$. We define $\tilde v(z)=v(z_0+r_0z)$ and $\tilde A(z)=r_0A(z_0+r_0z)$. Then 
\[
\Delta \tilde v+\nabla\cdot(\tilde A\tilde v)=0\quad\text{in}\quad B_8
\]
and
\[
\int_{B_8}|\tilde A|^2dz=\int_{B_{8r_0}(z_0)}|A|^2dz\le C_0.
\]
Hence, we have that $\tilde v(z)=0$ in $B_1$, namely, $v=0$ in $B_{r_0}(z_0)$. Using similar arguments as in the proof of Corollary~\ref{sucp}, we then conclude that $v$ is identically zero in $\Omega$.\eproof


\begin{thebibliography}{9999999}

%\bibitem[Al66]{al66}
%A. D. Aleksandrov, \emph{Majorization of solutions of second-order linear elliptic equations}, Vestnik Leningrad Univ., \textbf{21} (1966), 5-25 (in Russian). English transl. in Amer. Math. Soc. Transl. (2), \textbf{68} (1968), 120-143.

\bibitem[BK05]{bk05}
J. Bourgain and C. Kenig, \emph{On localization in the Anderson-Bernoulli model in higher dimensions}, Invent. Math., \textbf{161} (2005), 389-426.

%\bibitem[CS99]{cs99}
%J. Cruz-Sampedro, \emph{Unique continuation at infinity of solutions to Schr\"odinger equations with complex-valued potentials}, Proc. Edinburgh Math. Soc., (2) \textbf{42} (1999), 143-153.

\bibitem[Da12]{da12}
B. Davey, \emph{Some quantitative unique continuation results for eigenfunctions of the magnetic Schršdinger operator}, Comm. PDE, to appear.

%\bibitem[FS72]{fs72}
%C. Fefferman and E. Stein, \emph{$H^p$ spaces of several variables}, Acta Mathematica, \textbf{129} (1972), 137-193.

%\bibitem[GT83]{gt83}
%D. Gilbarg and N. S. Trudinger, \emph{Elliptic Partial Differential Equations of Second Order}, Springer, New York (1983). 

%\bibitem[JN61]{jn61}
%F. John and L. Nirenberg, \emph{On functions of bounded mean oscillation}, Comm. Pure Appl. Math., {\bf 14} (1961), 415-426.

%\bibitem[Ke93]{ke93}
%C. Kenig, \emph{Potential theory of non-divergence for elliptic equation}, "Dirichlet Forms", Lecture Notes in Math., {\bf 1562} (1993), 89-128.

\bibitem[KN00]{kn00}
C. Kenig and N.  Nadirashvili, \emph{A counterexample in unique continuation}, Math. Res. Lett., \textbf{7} (2000), 625-630.

\bibitem[KLW14]{klw14}
C. Kenig, L. Silvestre, and J.N. Wang, \emph{On Landis' conjecture in the plane}, preprint, arXiv:1404.2496v2 [math.AP]. 

\bibitem[Ki89]{ki89}
Y. M. Kim, \emph{Unique continuation theorems for the Dirac operator and the Laplace operator}, Doctoral thesis, MIT, 1989. 

\bibitem[KT02]{kt02}
H. Koch and D. Tataru, \emph{Sharp counterexamples in unique continuation for second order elliptic equations}, J. Reine Angew. Math., \textbf{542} (2002), 133-146.

\bibitem[KL88]{kl88}
V. A. Kondratiev and E. M. Landis,  \emph{Qualitative properties of the solutions of a second- order nonlinear equation}, Encyclopedia of Math. Sci. 32 (Partial Differential equations III), Springer-Verlag, Berlin (1988).

\bibitem[LW13]{lw13}
C. L. Lin and J. N. Wang, \emph{Quantitative uniqueness estimates for the general second order elliptic equations},  J. Funct. Anal., \textbf{266} (2014), 5108-5125.

\bibitem[Ma02]{ma02}
N. Mandache, \emph{A counterexample to unique continuation in dimension two}, Comm. Anal. Geom., \textbf{10} (2002), 1-10.

%\bibitem[Sa10]{sa10}
%M.V. Safonov, \emph{Non-divergence elliptic equations of second order with unbounded drift}, Nonlinear %Partial Differential Equations and Related Topics, AMS Transl. (2), \textbf{229} (2010), 211-232. 
\bibitem[St72]{st72}
R. S. Strichartz, \emph{A note on Trudinger's extension of Sobolev's inequality}, Indiana Univ. Math. J., \textbf{21} (1972), 841-842.

\bibitem[Tr67]{tr67}
N.S. Trudinger, \emph{On imbeddings into Orlicz spaces and some applications}, J. Math. Mech.,  \textbf{17} (1967), 473-483.

\bibitem[Wo90]{wo90}
T. Wolff, \emph{Unique continuation for $|\Delta u|\le V|\nabla u|$ and related problems}, Revista Math. Iberoamericana, \textbf{6} (1990), 155-200.

\bibitem[Wo94]{wo94}
T. Wolff, \emph{A counterexample in a unique continuation problem}, Comm. Anal. Geom., \textbf{2} (1994), 79-102. 

\bibitem[Wo93]{wo93}
T. Wolff, \emph{Recent work on sharp estimates in second order elliptic unique continuation problems}, The Journal of Geometric Analysis, \textbf{3} (1993), 621-650.

\bibitem[Ve62]{ve62}
I. N. Vekua, \emph{Generalized Analytic Functions}, Pergamon Press, London, 1962.

\end{thebibliography}
\end{document}